
%
%
%
%
%
%
%
%
\input amstex
\documentstyle{amsppt}
\magnification \magstep1
\parskip 11pt
\parindent .3in
\pagewidth{5.2in}
\pageheight{7.2 in}

%
\def \bl{\vskip 11pt}

\def \ni{\noindent}
\def \P{\bold{P}}
\def \Q{\bold{Q}}

\def \C{\bold{C}}
\def \Sig{\Sigma}
\def \mult{\text{mult}}
\def  \alb{\text{alb}}
\def \Alb{\text{Alb}}
\def \Pic{\text{Pic}^0}
\def \w{\omega}
\def \O{\Cal{O}}
\def \lra{\longrightarrow}
%
%

\centerline{\bf  SINGULARITIES OF THETA DIVISORS,}
\vskip 3pt
\centerline{\bf AND THE BIRATIONAL GEOMETRY OF IRREGULAR
VARIETIES}

\vskip 10pt

\centerline{ Lawrence EIN\footnote{Partially supported
by N.S.F. Grant DMS  93-02512   }}

\centerline{ Robert LAZARSFELD\footnote{Partially
supported by N.S.F. Grant  DMS  94-00815}}

\vskip 10pt
\ni {\bf Introduction}

The purpose of this paper is to show how the generic
vanishing theorems of
\cite{GL1} and \cite{GL2} can be used to settle a number
of  questions and conjectures raised  in \cite{Kol3},
Chapter 17, concerning the geometry of irregular complex
projective varieties. Specifically, we focus on three
sorts of results. First, we establish  a well known
conjecture characterizing principally polarized abelian
varieties whose theta divisors are singular in
codimension  one. Secondly, we study the holomorphic
Euler characteristic of varieties of general type having
maximal Albanese dimension: we verify
 a conjecture of Koll\'ar for subvarieties of abelian
varieties, but show that it fails in general. Finally,
we  give a surprisingly simple new proof of  a
fundamental theorem of Kawamata
\cite{Ka} on the Albanese mapping of varieties of
Kodaira dimension zero.

	 Turning to a more detailed description, we start  with
the singularities of theta divisors. Let
$A$ be an  abelian variety of dimension $g \ge 2$, and
let $\Theta
\subset A$ be a principal polarization on $A$, i.e. an
ample divisor such that
$h^0(A, \O_A(\Theta)) = 1$. Ever since the work
\cite{AM} of Andreotti and Mayer on the Schottky
problem, there has been  interest  in understanding what
sort of singularities $\Theta$ can have. A well-known
theorem of Kempf
\cite{K} states that if $A$ is a  Jacobian, then $\Theta$
has only rational singularities. For an arbitrary
principally polarized abelian variety $(A , \Theta)$,
Arbarello and DeConcini \cite{AD}  conjectured that if
dim Sing$(\Theta) = g - 2$ then  $(A , \Theta)$ splits
as a non-trivial product, i.e. that there
 exist principally polarized abelian varieties
$(A_1 , \Theta_1)$ and $(A_2 , \Theta_2)$ such that
$$A = A_1 \times A_2 \ \ \text{and} \ \ \Theta = pr_1^*
\Theta_1 + pr_2^* \Theta_2.$$  Koll\'ar \cite{Kol3,
Theorem 17.13}  recently put these matters into a new
perspective by proving that the pair $(A ,
\Theta)$ is log canonical. Denoting by
$\Sigma_k(\Theta)$ the multiplicity locus
$$ \Sigma_k(\Theta) = \left \{  x \in A \mid
\mult_x(\Theta) \ge k
\right \},
$$
 Koll\'ar's theorem implies in particular that
$$
\text{ every component of $\Sig_k(\Theta)$ has
codimension $\ge k$ in
$A$.}
$$  For example, taking $k = g+1$ it follows that
$\Theta$ can have no points of multiplicity $> g$.
Using a very pleasant enumerative argument, Smith and
Varley \cite{SV} subsequently established that if
$\Theta$ contains a $g$-fold point, then $(A ,
\Theta)$ splits as a product of $g$ elliptic curves (see
\cite{Nak} for a somewhat different approach).

	Our first result shows that the conclusion  of Kempf's
theorem holds quite generally, and that in fact any
example on the boundary of Koll\'ar's theorem is split:
\proclaim{\bf Theorem 1} If $\Theta \subset A$ is an
irreducible theta divisor, then $\Theta$ is normal and
has only rational singularities.
\endproclaim
\proclaim{\bf Corollary 2}  If $(A, \Theta)$ is any
principally polarized abelian variety, and if $k
\ge 2$, then $\Sigma_k(\Theta)$ contains an irreducible
component of codimension
$ \ k $ in $A$ if and only if
$(A , \Theta)$ splits as a
 $k$-fold product of p.p.a.v.'s.
\endproclaim
\ni When $k = g$ we recover the theorem of Smith and
Varley; the case
$k = 2$ gives the  conjecture of Arbarello and DeConcini.

The proof of Theorem 1 is surprisingly quick. In brief,
let $X \lra
\Theta$ be a resolution of singularities. By applying
the generic vanishing theorems  on $X$, and arguing with
some Nadel-type adjoint ideals on $A$, one reduces to
showing that $\chi(X,
\omega_X) > 0$. But $X$ is of general type, and the
inequality in question emerges as a special case of a
conjecture of Koll\'ar, which we discuss next.

Consider then a smooth projective variety $X$ of
dimension $n$, and assume that the Albanese mapping
$$\alb_X : X \lra \Alb(X)$$ is generically finite, or in
other words that $X$ has maximal Albanese dimension. A
result of
\cite{GL1} asserts that under these circumstances
$\chi(X, \omega_X)
\ge 0$. If $X$ is birationally the product of a torus
and some other variety then of course $\chi(X,\omega_X)
= 0$. But Koll\'ar conjectured \cite{Kol2, (17.9)} that
if $X$ is of general type, then
$\chi(X, \omega_X) > 0$.  Our second result shows that
this is true  if
$X$ is birationally a subvariety of
$\Alb(X)$:
\proclaim{\bf Theorem 3}\footnote{This proof of this
result builds on discussions some years ago with M.
Green, and it should be considered at least partially
joint work with him.} Let
$X$ be a smooth projective variety of maximal Albanese
dimension, and suppose that
$\chi(X,\omega_X) = 0$. Then the Albanese image {\rm
$$\alb_X(X)
\subseteq \Alb(X) $$}of
$X$ is fibred by tori. In particular, if {\rm
$\alb_X : X \lra \Alb(X)$} is birational onto its image
then $X$ is not of general type.
\endproclaim
\ni This is more than enough to give the inequality
required for Theorem  1.\footnote{In the interests of
truth in advertising,  we note that for Theorem 1  one
only needs Koll\'ar's conjecture when $X$ is birational
to an irreducible theta divisor, and this case is
covered by an old result of Kawamata and Viehweg
\cite{KV}. So in fact
 Theorem 1 can be read independently of the rest of the
paper. However we naturally prefer to see it as part of
a broader picture.} For smooth subvarieties of abelian
varieties, an equivalent statement was etablished
independently by Qi Zhang \cite{Z}.  We complete the
picture by showing that  Koll\'ar's conjecture fails  in
general. Our example is a threefold whose Albanese
mapping is a branched covering with a rather degenerate
branch divisor.

Finally, we turn to varieties of Kodaira dimension zero.
Let $X$ be a smooth projective variety of dimension $n$.
Kawamata \cite{Ka} showed that if $\kappa(X) = 0$, then
the Albanese mapping $\alb_X : X \lra
\Alb(X)$ is surjective. By a standard covering argument,
it is enough to prove this  assuming that $P_1(X) \ne
0$, where as usual
$P_{m}(X) = h^0(X, \omega_X^{\otimes m})$ denotes the
$m^{\text{th}}$ plurigenus  of $X$.  Kawamata's result
is therefore a consequence of
\proclaim{Theorem 4} If $P_1(X) = P_2(X) = 1$, then the
Albanese mapping of $X$ is surjective.\endproclaim
\ni Several other effective versions of Kawamata's
theorem were previously given by Koll\'ar (\cite{Kol1},
\cite{Mori}, \cite{Kol2},
\cite{Kol3}), the strongest of which states that
$\alb_X$ is surjective as soon as
$P_3(X) = 1$. Koll\'ar also asked for analogous results
involving
$P_2$.

However the main interest of Theorem 4 derives not so
much from any numerical improvements as from the the
surprising simplicity and transparency of its  proof.
Kawamata's approach involved some rather difficult
positivity results for  direct images of dualizing
sheaves, which were gradually replaced in Koll\'ar's
work by  subtle arguments with vanishing theorems. By
contrast, granting the general results of \cite{GL2},
Theorem 4 requires only a few  lines.  A pleasant
geometric argument also recovers the more precise
statement from
\cite{Ka} that if $\kappa(X) = 0$, then the fibres of
$\alb_X : X
\lra \Alb(X)$ are connected, as well as one of Koll\'ar's
characterizations of abelian varieties. We hope that
some of these ideas may find other applications in the
future.

The paper is organized as follows. In section 1 we
review for the convenience of the reader the results we
use from
\cite{GL1} and \cite{GL2}, and we prove Theorem 3. The
applications to varieties of Kodaira dimension zero
occupy \S 2, while
\S3 is devoted to  theta divisors.
We also include in \S 3 an extension of Koll\'ar's
theorem to singularities of pluri-theta divisors, as
proposed in \cite{Kol3}, Problem 17.15. We reiterate that
\S 3 may be read independently of the rest of the paper.

We have profitted from discussions with E. Arbarello, M.
Green, J. Koll\'ar, R. Smith, R. Varley and J. Wahl. In
particular, the statement of Theorem 1 was suggested by
Koll\'ar and Wahl, and as noted above discussions with
Green played a  substantial role in the proof of Theorem
3.

\vskip 10pt
\ni{\bf \S 0. Notation and Conventions}

\ni (0.1).  We work thoughout over the complex numbers
$\C$.

\ni (0.2).   Given a smooth variety or complex manifold
$X$ of dimension $n$, we generally denote by
$\w_X = \Omega^n_X$ the canonical line bundle of $X$. An
exception occurs in our discussion of adjoint ideals at
the beginning of \S 3, where in accordance with the
standard notation of higher dimensional geometry, we use
$K_X$ to denote (the linear equivalence class of) a
canonical divisor on $X$.

\ni (0.3).  If $D$ and $E$ are divisors on a variety or
complex manifold $X$, we write $D \prec E$ to indicate
that $E - D$ is effective.

\vskip 15pt

\ni {\bf \S 1.  Positivity of Holomorphic Euler
Characteristics}

We start by recalling the material from \cite{GL1} and
 \cite{GL2} that will be needed here and in \S2.  Then
we give the proof of Theorem 3  and the counter-examples
to the general case of Koll\'ar's conjecture.
\vskip 5pt

\ni {\it Review of Generic Vanishing Theorems}

Let $X$ be a compact K\"ahler manifold of dimension $n$,
and as usual let $\Pic(X)$ be the complex torus
parametrizing topologically trivial line bundles on $X$.
Given a point $y \in \Pic(X)$, we
 denote by $P_y$ the corresponding topologically trivial
bundle on
$X$. For $0 \le i \le n$  consider the analytic
subvarieties of
$V_i(X) \subset \Pic(X)$ defined  by:
$$
\aligned V_i(X) &= \left \{ y \in \Pic(X) \mid H^i(X,
\w_X \otimes P_y) \ne 0
\right \} \\
 &= \left \{ y \in \Pic(X) \mid H^{n-i}(X , P_y^*) \ne 0
\right
\}.
\endaligned $$ (For subsequent geometric arguments the
first description is preferable, but for the present
discussion the dual interpretation is easiest.)  Let
$y
\in V_i(X)$ be any point, and let
$$0 \ne v \in T_y  \ \Pic(X) = H^1(X, \O_X)$$ be a
non-zero tangent vector to $\Pic(X)$ at $y$. One of the
main themes of \cite{GL1} is that the first order
deformation theory of the groups $H^{n-i}(X, P_y^*)$ is
governed by the {\it derivative complex}:
$$  H^{n-i-1}(X , P_y^*) @> \cup v >> H^{n-i}(X , P_y^*)
@> \cup v >> H^{n-i+1}(X , P_y^*). \tag 1.1
$$ Roughly speaking, if $y \in V_i(X)$ is a sufficiently
general point, then
$v \in T_y \Pic(X)$ is tangent to
$V_i(X)$ if and only if the two maps in (1.1) vanish,
whereas  if (1.1) is exact, then all the cohomology in
$H^{n-i}(X, P_y^*)$ vanishes to first order in the
direction of $v$.  The principal result of \cite{GL2} is
that there are no  higher obstructions to deforming the
cohomology of topologically trivial line bundles, so
that a first order statement is equivalent to a global
one.

More precisely, one has the following
\proclaim{Theorem  1.2} {\rm (\cite{GL1}, \cite{GL2})}.
Fix  any irreducible component
$$S \subset V_i(X),$$ and let $y \in S$ be a general
point, i.e. a smooth point of $V_i(X)$ at which the
function $h^{n-i}(X, P_y^*)$ assumes its generic value
on $S$. Then:

\ni {\rm (1.2.1).} $S$ is a translate of a subtorus of
$\Pic(X)$.

\ni {\rm (1.2.2).}  $ \text{ \rm codim}_{\Pic(X)} S \ge
i -
\left( \dim X - \dim alb_X (X) \right )$.

 \ni {\rm (1.2.3).} If $0 \ne v \in H^1(X, \O_X)$ is
tangent to $S$, then the maps in {\rm (1.1)} vanish,
whereas if $v$ is not tangent to
$S$ then {\rm (1.1)} is exact.
\endproclaim
\demo{Proof} The first assertion is Theorem 0.1 of
\cite{GL2}, and the second is \cite{GL1}, Theorem 1.
For (1.2.3), let $\Delta_v(y)
\subset
\Pic(X)$ be a neighborhood of $y$ in the ``straight
line" in $\Pic(X)$ through $y$ determined by
$v$, i.e. the image of a small disk under the
exponential mapping
$\text{exp}_v(y) :
\C \lra \Pic(X)$ based at $y$ in the direction $v$. Thus
by (1.2.1), $\Delta_v(y) \subset S$ if
$v$ is tangent to $S$. Corollary 3.3 of \cite{GL2}
asserts that for
$t \in \Delta_v(y)$ in  punctured neighborhood of $y$:
$$ h^{n - i}(X, P^*_t) = \ \text{ dimension of homology
of (1.1)} .\tag 1.2.4
$$ Now $h^{n-i}(X, P_y^*)$ assumes its generic value at
$y$, and hence if
$\Delta_v(y) \subset S$ then $$h^{n-i}(X, P_t^*) =
h^{n-i}(X, P^*_y)$$ for all $t \in \Delta_v(y)$. It then
follows from (1.2.4) that the maps in (1.1) must vanish.
Similarly, if $\Delta_v(y) \not \subseteq S$, then the
left-hand side of (1.2.4) vanishes for generic $t$
thanks to the fact that $S$ is an irreducible component
of $V_i(X)$, and consequently (1.1) is exact. \qed

Still following \cite{GL1} and \cite{GL2}, Theorem 1.2
becomes particularly useful if it is restated via Hodge
duality. After fixing a K\"ahler metric on $X$, Hodge
theory gives conjugate linear isomorphisms:
$$H^{n-i}(X, P_y^*) \cong H^0(X, \Omega^{n-i}_X \otimes
P_y). \tag 1.3 $$  Similarly, if we represent $v \in
H^1(X, \O_X)$ by a harmonic
$(0,1)$-form, then  its conjugate is a holomorphic
one-form $\eta =
\bar v \in H^0(X,\Omega^1_X)$, and the conjugate of
(1.1) is:
$$H^0(\Omega^{n-i-1}_X \otimes P_y) @> \wedge \eta >>
H^0(\Omega^{n-i}_X \otimes P_y) @> \wedge \eta >>
H^0(\Omega^{n-i+1}_X \otimes P_y). \tag 1.4
$$
\proclaim{\bf Corollary 1.5} {\rm (1.5.1).} Keep
notation  and assumptions as in Theorem 1.2, and let
$\eta = \bar v$. Then
 $v$ is tangent to $S$ if and only if the maps in (1.4)
vanish, and otherwise (1.4) is exact.

\ni{\rm (1.5.2).} Assume that  $H^0(X, \omega_X) \ne 0$,
so that that the origin $0 = [\O_X] \in \Pic(X)$ lies
in  $V_0(X)$.  Then it  is an isolated point of $V_0(X)$
if and only if for every non-zero
$\eta \in H^0(X,\Omega^1_X)$, the map
$$ H^0(X , \Omega^{n-1}_X) @> \wedge \eta >> H^0(X,
\Omega^n_X)
$$ determined by wedging with $\eta$ is surjective.
\endproclaim
\demo{Proof} The first assertion is merely a restatement
of (1.2.1), and (1.5.2) follows by taking $i = 0$. \qed
\enddemo

\ni{\bf Remark 1.6}. The following slight generalization
will be useful. Let $$a : X \lra A$$ be a holomorphic
mapping from $X$ onto some complex torus $A$, and define
$$V_i(X)_A = \left \{ y \in \Pic(A) \mid H^{n-i}(X, a^*
P_y) \ne 0
\right \}.$$ Then the evident analogues of (1.2) and
(1.5) hold for these loci, where one works in $\Pic(A)$
instead of $\Pic(X)$ (so that the right-hand side of
(1.2.2) involves dim $a(X)$ rather than
$\dim \alb_X(X)$) and where the holomorphic one-forms
that occur in (1.4) and (1.5) are the pull-backs of the
flat one-forms on $A$. This is not stated explicitly in
\cite{GL1} and \cite{GL2}, but it is the natural context
in which the arguments there work.

\vskip 10pt
\ni{\it Holomorphic Euler Characteristics}

We assume for the remainder of this section that $X$ is
a compact K\"ahler manifold of dimension $n$ whose
Albanese mapping
$$a = \alb_X : X  \lra \Alb(X) = A$$ generically finite
onto its image. It follows from \cite{GL1, Theorem 1},
or (1.2.2) above, that then
$$h^0(X, \omega_X \otimes P^*_y) = \chi(X, \omega_X)
\tag 1.7$$ for general $y \in \Pic(X)$.  Our goal is to
show that if
$\chi(X, \omega_X) = 0$, then the Albanese image of $X$
is ruled by tori.

We start with the following  useful remark due to M.
Green:
\proclaim{Lemma 1.8}  One has inclusions:
$$\Pic(X) \supseteq V_0(X) \supseteq V_1(X) \supseteq
\dots \supseteq V_{n}(X) = \{ \O_X \}. $$ \endproclaim
\demo{Proof} Suppose that $y \in V_{i}(X)$ ($i > 0$), so
that $H^{n-i}(X, P_y^*) \ne 0$. In view of (1.3), there
exists a non-zero form $0 \ne  \alpha \in  H^0(X,
\Omega^{n-i}_X \otimes P_y)$. It is enough to show
that $\eta
 \wedge \alpha \ne 0 \in H^0(X,\Omega^{n-i+1}_X \otimes
P_y)$ for general $\eta \in H^0(X, \Omega_X^1).$ But
this follows as in
\cite{GL1, end of proof of Theorem 2.10}.  In
brief, fix a general point $x
\in X$ at which $\alpha(x) \ne 0$. Since $X$ has maximal
Albanese dimension, if we have choosen $x$ sufficiently
generally, we can find holomorphic one-forms
$\eta_1, \dots ,\eta_n \in H^0(X, \Omega^1_X)$ such that
the  $
\eta_i(x) $ form a basis of the holomorphic cotangent
space $T^*_x X$. But then it is immediate that
$\alpha(x) \wedge \eta_j(x) \ne 0$ for some $j \in [1,
n]$. \qed

We now turn to the demonstration of Theorem 3. As noted
in the Introduction, the argument that follows builds
substantially on discussions with M. Green some years
ago.
\demo{Proof of Theorem 3}  Assume that $\chi(X,
\omega_X) = 0$. Then by (1.7), $V_0(X)$ is a proper
subvariety of $\Pic(X)$. Fix an irreducible component $S
\subset V_0(X)$ and  a general point $y \in S$,  and put
$$k = \text{codim}_{\Pic(X)} S. $$ Note that it follows
from (1.2.2) that $S$ cannot be contained in
$V_{j}(X)$ for  $j > k$. By contrast,  if $S \subseteq
V_j(X)$ for some
$j \le k$, then in fact $S$ is an irreducible component
of $V_j(X)$ thanks to the previous lemma. Therefore,  by
taking $y \in S$ sufficiently generally, we may suppose
that
$H^{n-j}(X, P^*_y) = 0$ for $j > k$, and that if
$H^{n-j}(X, P^*_y) \ne 0$ for $j \le  k$, then the
conclusions of Theorem 1.2 hold with $i = j$ at $y$. In
particular, it follows  from (1.2.3) that if $0
\ne v \in H^1(X,\O_X) = T_y \Pic(X)$ is not tangent to
$S$ at $y$, then the  derivative complex
$$0 @>>> H^{n-k}(P_y^*) @> \cup v >> H^{n-k+1}(P^*_y) @>
\cup v >>
\dots @> \cup v >> H^n(P^*_y) @>>> 0 \tag 1.9$$ is
(everywhere) exact.

We claim that $$H^{n-k}(X, P^*_y) \ne 0, \tag 1.10$$
i.e. that $S$ is actually a component of $V_k(X)$ (and hence
also that $k \le n$). In fact, let $$V
\subset H^1(X, \O_X) = T_y  \ \Pic(X)$$ be a
$k$-dimensional subspace complementary to
$T_y S \subset T_y \Pic(X)$ (so that $V$ represents the
normal space to
$S$ at $y$.) Thus (1.9) is exact for each $0 \ne v \in
V$. Set $\P =
\P(V^*)$, so that
$\P$ is a projective space of dimension $k-1$. Then we
may assemble the derivative complexes (1.9) determined
by $0 \ne v \in V$ into a complex $K$ of vector bundles
on $\P$:
$$0 \rightarrow H^{n-k}(P^*_y)\otimes \O_\P(-k)
\rightarrow H^{n-k+1}(P^*_y) \otimes \O_\P(-k + 1)
\rightarrow \dots
\rightarrow H^n(P_y^*) \otimes \O_\P
\rightarrow 0.
$$ The fact that each of the point-wise complexes (1.9)
is exact implies that $K$ is exact as a complex of
sheaves on $\P$.  It then follows by chasing through
$K$ and taking cohomology on $\P$ that
$$0 \ne H^n(X,  P^*_y) \cong H^{n-k}(X, P^*_y),$$  and
(1.10) is established. [Compare
\cite{Mori}, Proof of (3.3.2).]

We now argue as in \cite{GL2}, \S 4. Recalling that $S$
is a (translate of) a subtorus of $\Pic(X)$, let $C =
S^*$ be the dual torus. Since $A = \Alb(X)$ is the dual
of $\Pic(X)$, the inclusion $S \hookrightarrow \Pic(X)$
determines a quotient map
$$\pi : A = \Alb(X) \lra C$$ whose fibres are translates
of the
$k$-dimensional connected subtorus $$B =_{\text{def}}
\ker (\pi).$$  Let $Y = \alb_X(X) \subset A$ be the
Albanese image of
$X$, and let $$
\gathered g  : X @> \alb_X >> Y \subset A  @> \pi >> C \\
h : Y \subset A  @> \pi >> C
\endgathered \tag 1.11
$$ denote the indicated compositions. We claim that
$$\dim g(X) \le n - k. \tag 1.12$$ Grant this for a
moment. Since $a : X \lra Y$ is generically finite and
surjective, it then follows that all the fibres of
$h : Y \lra h(Y) \subset C$ have dimension $\ge k$. But
these fibres are contained in translates of the
$k$-dimensional torus $B$. In other words,  the fibres of
$Y \lra h(Y)$ fill up the fibres of $A \lra C$ over
$h(Y)$. Therefore $Y$ is ruled by tori, as was to be
shown.

It remains only to prove (1.12). But this is in fact
established in
\cite{GL2, p. 92}. We summarize the argument for the
convenience of the reader. Let
$$v_1, \dots , v_{q-k} \in H^1(X, \O_X) = T_y S$$ be a
basis for the tangent space to $S$ at $y$, where
 $q = \dim \Pic(X)$. As in the previous subsection, let
$\eta_i = \bar v_i \in H^0(X, \Omega_X^1)$ be conjugate
holomorphic one-forms. The map $g : X \lra C$ arises by
integrating the $\eta_i$, and consequently for a general
point $x \in X$:
$$\dim g(X)  = \dim \text{ span} \left \{ \eta_1(x),
\dots,
\eta_{q-k}(x) \right \} \subset T^*_x X.$$ Since
 $\eta_i$ is the conjugate of a tangent vector to $S$,
it follows from (1.5.1) that each of the maps
$$H^0(X, \Omega^{n-k}_X \otimes P_y) @> \wedge \eta_i >>
H^0(X,
\Omega^{n-k+1}_X \otimes P_y),
\tag *$$ vanishes.  But the group on the left in (*) is
non-zero thanks to (1.10). An elementary pointwise
calculation (cf.
\cite{GL2}, Lemma 4.1) shows that the space of one-forms
that annihilates a non-zero $n-k$ form has dimension
$\le n - k$, and hence (1.12) follows. \qed

We conclude this section by giving an example to show
that Koll\'ar's conjecture \cite{Kol3, (17.9)}  fails in
general.

\ni{\bf Example 1.13}. We exhibit a smooth variety $X$
of general type having maximal Albanese dimension such
that $\chi(X, \omega_X)= 0$. Let $E$ be an elliptic
curve, and $p : C \lra E$ a double covering of $E$ by a
curve $C$ of genus $\ge 2$. Denote by $\iota : C
\lra C$ the corresponding involution of $C$. Let
$A = E \times E \times E$,
$V = C \times C \times C$, and consider the involution
$\tau = \iota
\times \iota \times \iota$ of $V$. Set $Y = V / \{ 1,
\tau \}$. The spaces in question fit into a tower of
Galois covers:
$$V @>f>> Y @>g>> A$$ of degrees 2 and 4 respectively.
Observe that $Y$ is smooth except at finitely many
points at which it is locally analytically isomorphic to
a Veronese cone, i.e. the quotient of $\C^3$ by
multiplication by $-1$. In particular $Y$ has only
terminal and hence rational singularities. The map $f :
V \lra Y$ being \'etale off the finite many singular
points of $Y$, we see that
$Y$ is of general type, and in fact minimal of global
index two. Let
$$h: X \lra Y$$ be a resolution of singularities.
Clearly $X$ has maximal Albanese dimension, and we claim
that $\chi(X, \omega_X) = 0$. In order to verify this,
it suffices in view of (1.7) to show that
$H^3(X, h^*g^* P) = 0$ for a general topologically
trivial line bundle $P \in \Pic(A)$. To this end, start
with non-trivial line bundles $P_1 , P_2, P_3 \in
\Pic(E)$, and take $P$ to be their exterior product.
Recalling that $Y$ has only rational singularities, one
calculates:
$$H^3(X, h^*g^* P) = H^3(Y, g^*P) = H^3(V,
f^*g^*P)^{\{1, \tau \}} = 0,$$ as required.
[Alternatively, as Koll\'ar notes, one can compute $(g \circ
h)_*(\w_X)$.]

\vskip 10pt

\ni {\bf \S 2.  Varieties of Kodaira Dimension Zero}

We give in this section the applications to varieties of
Kodaira dimension zero. Denote by $X$  a smooth projective
variety of dimension
$n$,  and  write $P_m(X)$ for the
$m^{\text{th}}$ plurigenus of $X$, i.e. $P_m(X) = h^0(X,
\omega_X^{\otimes m})$. As in \S 1 we consider  the loci:
$$V_i(X) = \left \{ y \in \Pic(X) \mid H^i(X, \omega_X
\otimes P_y)
\ne 0 \right \}.$$
 We emphasize that we do not  assume here that
the Albanese mapping of $X$ is generically finite.

Theorem 4
from the Introduction  is an
immediate consequence of two elementary general propositions:
\proclaim{\bf Proposition 2.1}  If $P_1(X) = P_2(X) =
1$, then the origin is an isolated point of $V_0(X)$.
\endproclaim
\proclaim{\bf Proposition  2.2}  If the origin is an
isolated point of
$V_0(X)$, then the Albanese mapping
$\alb_X : X \lra \Alb(X)$ is surjective. \endproclaim

\demo{Proof of Proposition 2.1} Since $P_1(X) \ne 0$, the
origin $\O_X$ lies in $V_0(X)$. Suppose that it is not
an isolated point. Then  by (1.2.1), $V_0(X)$ contains a
subgroup $S \subset V_0(X)$ of positive dimension. In
particular, if $y \in S$ then also
$-y \in S$, and therefore for each $y \in S$ the image
of
$$H^0(X,\omega_X \otimes P_y) \otimes H^0(X,\omega_X
\otimes P_y^*)
\lra H^0(X, \omega_X^{\otimes 2}) \tag *$$ is non-zero.
Since a given divisor in the linear series
$|\omega_X^{\otimes 2}|$ has only finitely many
irreducible components, it follows that as $y$ varies
over the positive dimensional torus
$S$, the image of (*) must vary as well. Therefore
$P_2(X) > 1$, a contradiction. [Compare the proof of
\cite{Kol3}, Theorem 17.10.]
\qed \enddemo

\demo{Proof of Proposition 2.2} Assume that the origin
is an isolated point of $V_0(X)$, but that $\alb_X$ is
not surjective.  Fix an arbitrary point $x \in X$. The
non-surjectivity of $\alb_X$ implies that  there exists a
non-zero holomorphic one-form
$$0 \ne \eta = \eta_x \in H^0(X, \Omega^1_X) \ \
\text{such that } \
\eta(x) = 0.$$ [Take $\eta$ to be the pull-back of a
flat one-form on $\Alb(X)$ lying in the kernel of the
coderivative $ T^*_{\alb(x)}\Alb(X)
\lra T^*_x X$ of $\alb_X$ at $x$.]  On the other hand,
since the origin is an isolated point of $V_0(X)$, it
follows from (1.5.2) that wedging with $\eta_x$ gives a
surjective map
$$H^0(X, \Omega^{n-1}_X) @> \wedge \eta_x >> H^0(X,
\Omega^n_X).$$ Therefore every section of
$\omega_X = \Omega^n_X$ vanishes at (the arbitrary
point) $x$, and hence $H^0(X, \omega_X) = 0$, a
contradiction. \qed \enddemo

\ni {\bf Remark 2.3}  For the convenience of the reader,
we recall the standard argument showing that Theorem 4
implies Kawamata's result that if $Y$ is a smooth
projective variety with $\kappa(Y) = 0$, then
$\alb_Y$ is surjective. In fact, by a lemma of Fujita
(cf.
\cite{Mori}, (4.1)), there exists a smooth projective
variety $X$ of Kodaira dimension zero, admitting a
generically finite surjective map $f : X \lra Y$, such
that $P_1(X)
\ne 0.$  Since $\kappa(X) = 0$, it follows that $P_1(X)
= P_2(X) = 1$. Therefore $\alb_X$ is surjective thanks
to Theorem 4, and this easily  implies that $\alb_Y$ is
surjective.

In the remainder of this section, we show how similar
ideas lead to new proofs of some other results
concerning varieties of Kodaira dimension zero. We start
with a theorem of Koll\'ar \cite{Kol1} giving a
birational characterization of abelian varieties.
Koll\'ar's statement was in turn an effective version of
a theorem of Kawamata-Viehweg \cite{KV} characterizing
abelian varieties birationally as projective manifolds
with $\kappa = 0$ and $q = n$. Some stronger results
appear in \cite{Kol2} and \cite{Kol3}, but it is not clear
whether one could  recover them by these techniques.

\proclaim{\bf Proposition 2.4} {\rm (Koll\'ar,
\cite{Kol1}).} Let $X$ be a smooth projective variety  of
dimension $n$, and assume that
$$P_1(X)  = P_4(X) = 1 \ \ \text{ and} \ \ \ q(X)
=_{\text{def}}h^1(X,\O_X) = n.$$ Then $X$ is
birational to an abelian variety. \endproclaim
\demo{Proof}  The beginning of the argument follows
Koll\'ar's proof. The Albanese mapping
$$a = \alb_X : X \lra \Alb(X) = A$$ is surjective by
Theorem 4, hence generically finite since $q = \dim A =
\dim X = n$. If $a$ is a birational isomorphism, we
are done. Otherwise, since the Albanese map does not
factor through any \'etale covers of $A$, it follows by
considering the Stein factorization of $a$ that the
ramification divisor of $a$ must contain at least one
component $\Delta$ which maps to a divisor $D \subset
A$. Replacing $X$ if necessary by a suitable blow-up, we
may assume that $\Delta$ is smooth.  Note that since the
ramification divisor $Ram(a)$ represents $\omega_X$, there is a
natural inclusion $\omega_X(\Delta) \hookrightarrow
\omega_X^{\otimes 2}$.

We now apply (1.2) and (2.2) to $\Delta$. Specifically,
since
$a(\Delta) = D$ is a divisor in $A$, we have
$P_1(\Delta) \ge  1$. Suppose first that $D$
spans $A$. Then evidently $\alb_\Delta$ is not
surjective.  Proposition 2.2 implies that $\O_\Delta$ is
not an isolated point of
$V_0(\Delta)$, which therefore contains a subtorus of
positive dimension. In fact, combining Remark 1.6
(applied to $\Delta$) with the proof of Proposition 2.2,
we see that there is a positive dimensional subgroup
$S_\Delta \subset \Pic(X)$ such that
$$H^0(\Delta, \omega_\Delta \otimes P_y ) \ne 0 \ \
\forall \ y \in S_\Delta. \tag 2.4.1$$

The hypotheses on the plurigenera of $X$ force $P_2(X) =
1$, and so Proposition 2.1 implies that $\O_X$ is an
isolated point of
$V_0(X)$. Therefore, thanks to Lemma 1.8, $\O_X$ is an
isolated point of all the $V_i(X)$, and in particular
$$H^0(X, \omega_X \otimes P_y) = H^1(X, \omega_X \otimes
P_y) = 0$$ for $y$ in  punctured neighborhood of $0$ in
$S_\Delta$. In view of (2.4.1), it  follows from the
exact sequence
$$ 0 \lra \omega_X \otimes P_y \lra \omega_X (\Delta)
\otimes P_y
\lra \omega_\Delta \otimes P_y \lra 0$$ that $H^0(X,
\omega_X(\Delta) \otimes P_y) \ne 0$ for every $y \in
S_\Delta$. Then just as in the proof of Proposition 2.1,
this implies that $$h^0(X, \omega_X^{\otimes 2}(2
\Delta)) \ge 2. $$  But $H^0(X, \omega_X^{\otimes 2}(2
\Delta)) \subset H^0(X, \omega_X ^{\otimes 4})$, and
hence $P_4(X) > 1$, a contradiction. It remains to treat
the possibility that $\Delta$ maps to a codimension one
subtorus $D \subset A$, but in this case it is enough to
take $S_\Delta$ to be the kernel of the natural map
$\Pic(X) = \Pic(A) \lra \Pic(\Delta)$. \qed \enddemo

We conclude this section by showing how similar ideas
lead to the more precise result from
\cite{Ka} that the Albanese mapping of a variety of
Kodaira dimension zero is a fibre space, i.e. has
connected fibres:
\proclaim{\bf Proposition 2.5} {\rm (Kawamata \cite{Ka}).  }
Let
$X$ be a smooth projective variety of Kodaira dimension
zero. Then the fibres of the Albanese mapping
$$a = \alb_X : X \lra \Alb(X) = A$$ are connected.
\endproclaim
\demo{Proof} There is no loss in generality in replacing
$X$ with a birationally equivalent model. So by
starting with the Stein factorization of $a$ and
resolving singularities and indeterminacies, we can
assume that  $a = \alb_X$ admits a factorization:
$$X @> g >> V @>b>> A,$$ where $g$ is surjective with
connected fibres, $V$ is smooth and projective, and $b$ is
generically finite.  We know already (e.g. by
Theorem 4 and Remark 2.3) that $a$ is surjective, and
hence $b$ is a generically finite covering. We assume by
way of contradition that
$b$ has degree $> 1$. Since the Albanese mapping does
not factor through any \'etale coverings of
$A$,  $b$ cannot be birationally \'etale, and
$V$ cannot be birational to an abelian variety.
Therefore $b$ has a non-trivial ramification divisor
$R$, and by Proposition 2.4, $P_4(V) > 1$.

We next reduce in effect to the situation $P_1 \ne 0$.
In fact, by Fujita's lemma \cite{Mori, (4.1)} there is
a smooth projective variety
$Y$, and a generically finite surjective covering
$\nu : Y \lra X$ such that $\kappa(Y) = \kappa(X) = 0$
and $P_1(Y) \ne 0$. Therefore $P_1(Y) = 1$, and we
denote by $K_Y$ the unique effective canonical divisor
on $Y$. We put
$f = g \circ \nu$, and consider the maps:
$$Y @> f >> V @>b >> A.$$ Let $\Delta \subset V$ be any
irreducible component of the ramification divisor $R =
Ram(b)$. We claim:
$$\text{ Any irreducible component $D$ of $f^*(\Delta)$
appears in
$K_Y$. }  \tag 2.5.1$$
Grant this for the time being. We
can write
$$f^*(R) = \sum a_i D_i \ \ \ \ (a_i > 0),
$$
where the $D_i$ are distinct irreducible divisors, and
$\sum D_i \prec K_Y$  by (2.5.1). Let $a = \max \{
a_i \}$. Then for any $m > 0$:
$$
\aligned f^*(mR) = \sum m a_i D_i &\prec \sum maD_i \\
        &\prec m a K_Y.
\endaligned \tag *
$$ On the other hand, since $b : V \lra A$ is a
generically finite covering of an abelian variety, one has
$R \equiv K_V$. Therefore $$P_4(V) = h^0(V,
\O_V(4R)) \ge 2,$$ and hence $h^0(Y,
\O_Y(4 a K_Y)) \ge 2$ thanks to (*). But $\kappa(Y) =
0$, so this is a contradiction.

It remains to prove (2.5.1). Let $v \in \Delta$ be any
point, and $y
\in f^{-1}(v)$. Since $b$ ramifies at $v$, as in the
proof of Propoisiton 2.2 there is a homolorphic one-form
$$0 \ne  \eta_v
\in H^0(V, \Omega^1_V) \ \ \text{such that  $\eta_v(v) =
0$.}$$
On the other hand, Proposition 2.1 and  the fact
that $\kappa(Y) = 0$ imply that $\O_Y$ is an isolated
point of $V_0(Y)$. It follows from (1.5.2) that the map
$$H^0(Y, \Omega^{n-1}_Y) @> \wedge (f^*\eta_v) >> H^0(Y,
\Omega^n_Y)$$ is surjective. But $f^*\eta_v (y) = 0$,
and hence $K_Y$ vanishes at
$y$. \qed \enddemo

\ni{\bf Remark 2.6.}  Propositions 2.1, 2.2 and 2.4 extend
with no difficulties to the setting of compact K\"ahler
manifolds. However it is not immediately obvious to us how
to arrange that the manifold $V$ in the proof of (2.5) be
K\"ahler.

\vskip 15pt
\newpage
\ni {\bf \S 3.  Singularities of Theta Divisors}

We start with some preliminaries on adjoint ideals, and
then give the applications to theta divisors. The reader
who wishes to read this section independently of the
rest of the paper is referred to Remark 3.4.

\vskip 5pt
\ni{\it Adjoint Ideals}
\def \J {\Cal J}

We wish to understand how the adjunction formula works
for a possibly singular divisor in a smooth variety. The
following Proposition generalizes various classical
constructions involving  conductor ideals. The result
is certainly at least implicitly known to the experts,
but we include a proof for the benefit of the reader.
\proclaim{\bf Proposition 3.1} Let $M$ be a smooth
variety, let $D \subset M$ be a
reduced effective divisor, and let
$f : X \lra D$ be any resolution of singularities. Then
there is an {\rm adjoint ideal}
$$\J  = \J_D \subset \O_M,$$ cosupported in the singular
locus of
$D$, which sits in an exact sequence:
$$ 0 @>>> \O_M(K_M) @> \cdot D >> \O_M(K_M + D) \otimes
\J @>>> f_*
\O_X(K_X)
\lra 0. \tag 3.1.1$$ Moreover, $\J = \O_M$ if and only
if $D$ is  normal and has only rational singularities.
\endproclaim
\demo{Proof}  Note to begin with that the sheaf $f_*
\w_X$ is independent of the choice of resolution, so we
are free to work with any convenient one.  Let
$g : Y \lra M$ be an embedded   resolution of $D$, and
let
$X \subset Y$ be the proper transform of $D$ (so $X$ is
a disjoint union of smooth divisors). Then we can write
$$K_Y + X = g^*(K_M + D) + P - N,$$ where $P, N$ and $X$
are effective divisors on $Y$, with no common components,
and every component of $P$ is $g$-exceptional.   The
adjoint ideal $\J = \J_D$ is defined to be $$\J = g_*
\O_Y(-N)
\subseteq g_*\O_Y = \O_M.$$ It follows from a lemma of
Fujita (cf. \cite{KMM, 1.3.2}) that
$g_* \O_P(P) = 0$, and consequently $g_* \O_P(P-N) = 0$.
We then have:
$$
\align g_* \O_Y(K_Y + X) &= \O_X(K_M + D) \otimes
g_*\O_Y(P-N) \\
                  &= \O_X(K_M + D) \otimes g_*\O_Y(-N) \\
                  &= \O_X(K_M + D) \otimes \J.
\endalign$$ Recalling that $g_* \O_Y(K_Y) = \O_M(K_M)$
and $R^1 g_* \O_Y(K_Y) = 0$,  (3.1.1) arises as the
pushforward under
$g$ of the exact sequence $$0 \lra \O_Y(K_Y) \lra
\O_Y(K_Y + X) \lra
\O_X(K_X) \lra 0.$$ It follows from (3.1.1) that $\J_D =
\O_M$ iff $f_* \w_X = \w_D$. Since $f$ factors through
the normalization of $D$, this can evidently  hold
only if $D$ is normal. And as $D$ is in any event
Cohen-Macaulay, the equality in question is then
equivalent to the condition that $D$ has at worst
rational singularities (cf. \cite{Kol4, (11.10)}).
\qed
\enddemo

\ni{\bf Remark 3.2.} As in the Introduction, consider the
multiplicity loci
$$\Sigma_k(D) = \left \{ x \in M \mid \mult_x(D) \ge k
\right \}.
\tag 3.2.2$$ We observe  that if
$\Sigma_k(D)$ has any components of codimension $\le k$
in $M$, for some $k \ge 2$, then the corresponding
adjoint ideal $\J = \J_D$ is non-trivial. In fact,
 construct the
embedded resolution $Y$ by first blowing up such a
component. Then one sees that the divisor $N$ appearing
in the proof of (3.1) is non-zero, and the assertion
follows. (Compare \cite{Kol3}, proof of Theorem 17.13.)

\vskip 15pt
\ni{\it Theta Divisors}

Let $(A, \Theta)$ be a principally polarized abelian
variety. We start with the following statement, which in
view of (3.1) is equivalent to Theorem 1 from the
Introduction.
\proclaim{\bf Theorem 3.3} Assume that $\Theta$ is
irreducible. Then the corresponding adjoint ideal is
trivial, i.e. $\J = \J_\Theta =
\O_A$.\endproclaim
\demo{Proof} We may assume that $g = \dim A \ge 2$, for
otherwise the statement is trivial. Let
$f:  X \lra \Theta$ be a resolution. Note first that
$X$ (i.e. $\Theta$) is of general type. Otherwise, by a
theorem of Ueno (cf.
\cite{Mori}, (3.7)) there would exist a non-trivial
quotient  $\pi : A \lra B$ of $A$, plus an effective
divisor $E
\subset B$ such that $\Theta \subset \pi^*(E)$. But this
is impossible since
$\Theta$ is ample. [We remark that Ueno's theorem is
the  essential point where irreducibility is used:  for
$\Theta$ reducible, the individual components won't be
of general type.]

The adjoint  exact sequence (3.1.1)  takes the form:
$$ 0 \lra \O_A \lra \O_A(\Theta) \otimes \J \lra f_*
\w_X \lra 0. \tag 3.3.1
$$ Now let $P \in \Pic(A) $ be a topologically trivial
line bundle on $A$, and twist (*) by $P$:
$$ 0 \lra P \lra \O_A(\Theta) \otimes P \otimes \J \lra
f_* \w_X \otimes P  \lra 0 . \tag 3.3.2
$$ Evidently $X$ has maximal Albanese dimension.
Bearing in mind Remark 1.6, it follows from
\cite{GL1, Theorem 1} (or (1.2.2) above) that
$H^i(X, \omega_X \otimes f^*P )  = 0$ for $i> 0 $ and
$P$ general. Therefore
$$
\aligned H^0(A, f_* \w_X \otimes P) &= H^0(X, \w_X
\otimes f^* P) \\ &= \chi(X, \w_X \otimes f^* P) \\ &=
\chi(X ,
\omega_X)
\endaligned
$$ for generic $P$. But by construction $X$ is
birational to its Albanese image, and hence
$  \chi(X, \w_X)
\ne 0$ thanks to Theorem 3. Thus $H^0(A, f_* \w_X \otimes
P) \ne 0$ for general $P$. Let $\alpha_P \in A$ be the
point corresponding to $P \in \Pic(A)$ under the
isomorphism $A \cong \Pic(A)$ determined by the
principal polarization $\O_A(\Theta)$. Then (3.3.2)
implies that
$$ H^0(A,  \O_A(\Theta + \alpha_P) \otimes \J ) =
H^0(A,  \O_A(\Theta)
\otimes P \otimes \J )    \ne 0 \tag 3.3.3$$ for general
$P \in \Pic(A)$ (and hence general $\alpha_P \in A$).

Suppose now that the theorem is false, so that  $\J   \ne
\O_A$. Then the corresponding zero-locus $$Z = \text{
Zeroes($\J$) }
\subset A $$ is non-empty. Since $h^0(A, \O_A(\Theta +
\alpha_P)) = 1$, it follows from (3.3.3) that
$Z \subseteq (\Theta + \alpha)$ for a general point
$\alpha \in  A$. But the translates of $\Theta$ don't
have any points in common, so this is impossible.
\qed
\enddemo

\ni{\bf Remark 3.4.} We indicate how to avoid the appeal
to Theorem 3 (and Remark 1.6) in the argument just
completed. Keeping notation as above, Theorem 3 was
invoked only in order to  show that
$\chi(X, \omega_X) \ne 0$. This can be circumvented by
noting  from (3.3.1) that
$H^0(X, f_* \omega_X) \subseteq H^1(A, \O_A)$, and hence
$p_g(X) \le g = \dim A$.   But then a theorem of Kawamata and
Viehweg \cite{KV} (cf. \cite{Mori, (3.11)}) implies that
$\chi(X, \omega_X) = 1$. In fact, Kawamata
and Viehweg prove  that the maps $H^0(A,
\Omega_A^i) \lra H^0(X, \Omega^i_X)$ are isomorphisms
for $i \le g-1$, and taking $i = 1$ it follows that
$\Pic(A)
\lra \Pic(X)$ is an isogeny. Hence one can also bypass
Remark 1.6: it is essentially the same to work with
general $P \in \Pic(X)$ (as in \cite{GL1}) or with
general $P \in \Pic(A)$ (as in the argument above).

\def \codim{\text{codim}}
 Next we give the
\demo{Proof of Corollary 2} Assume that there exists an
integer $k
\ge 2$ such that $\Sigma_k(\Theta)$ has an irreducible
component  of codimension $k$ in $A$. We wish to show
that then $(A, \Theta)$ is a  non-trivial $k$-fold
product of principally polarized abelian varieties.  It
follows in the first place from Remark 3.2 that
$\J_\Theta \ne \O_A$. Therefore $\Theta$ is reducible by
Theorem 3.3. The Decomposition Theorem (cf. \cite{LB,
(4.3.1)}) then implies that the
p.p.a.v. $(A, \Theta)$ splits as a non-trivial product.
Let
$$(A, \Theta) = (A_1, \Theta_1) \times \dots \times
(A_r, \Theta_r)$$ be the decomposition of $(A, \Theta)$
as a product of irreducible p.p.a.v.'s. Given any
irreducible subset $S \subset
\Sigma_k(\Theta)$, there exist integers $k_1, \dots ,
k_r \ge 0$, with $\sum k_i \ge k$, such that
$$
S \subseteq    \Sigma_{k_1}(\Theta_1)
\times \dots \times
\Sigma_{k_r}(\Theta_r)  .
$$  Suppose that $\codim_A S = k$. Since in any event
$\codim_{A_i}\Sigma_{k_i}(\Theta_i) \ge k_i$ by
Koll\'ar's theorem, it follows that
$$\aligned k =  \codim_A S &\ge
\codim_{A_1} \Sigma_{k_1}(\Theta_1) + \dots +
\codim_{A_r}
\Sigma_{k_r}(\Theta_r) \\ &\ge k_1 + \dots + k_r \\
 &\ge k.
\endaligned
\tag *$$  Therefore
$$\codim_{A_i}\Sigma_{k_i}(\Theta_i) = k_i  \ \ \text{
for all} \ \ i.$$ By induction, this implies that $k_i =
0 \text{ or } 1$ for all
$i$. But then it follows from (*) that $r \ge k$. \qed
\enddemo

\def \eps{\epsilon}

Finally, we present an extension of Koll\'ar's theorem
to the singularities of pluri-theta divisors, as
proposed in
\cite{Kol3, Problem 17.15}.
\proclaim{\bf Proposition 3.5} Let $(A, \Theta)$ be a
p.p.a.v., and for $m \ge 1$ fix any divisor $D \in |m
\Theta|$. Then the pair $(A,
\frac{1}{m}D)$ is log-canonical. In particular, every
component of
$\Sigma_{mk}(D)$ has codimension $\ge k$ in $A$.
\endproclaim
\demo{Proof} The point of the argument is to use
systematically all the available vanishings.
Specifically, for $0 < \epsilon \ll 1$ consider the
divisor $$E = E_\eps = \frac{1- m\eps}{m}D\equiv (1-\eps)
\Theta.$$  Thus $\Theta - E$ is an ample $\Q$-divisor,
and hence so is $P(\Theta - E)$ for any
$P \in \Pic(A)$. Denote by
$\J =\J_E \subset \O_A$ the multiplier ideal determined
by $E$ (cf.
\cite{De, (6.12)} or \cite{Kol4, (2.16)}). Then
Kawamata-Viehweg-Nadel vanishing implies that
$$H^i(A, \O_A(\Theta) \otimes P \otimes \J) = 0 \ \
\text{ for all }
\ i > 0 \ \text{ and } \ P \in \Pic(A).  \tag 3.5.1$$

Assuming the assertion of the Proposition false, we can
choose $\eps \ll 1$ such that $\J \ne \O_A$. Then $Z =
\text{ Zeroes}(\J) \ne \emptyset$, and as in the proof of
Theorem 3.3 it follows that $$H^0(A, \O_A(\Theta) \otimes
P \otimes \J) = 0 \tag 3.5.2$$ for general
$P \in \Pic(A)$. Combining (3.5.1) and (3.5.2), we find
that
$$\chi(A, \O_A(\Theta) \otimes P \otimes \J) = 0$$ for
general $P$. As Euler characteristics are deformation
invariants, this then holds  for {\it arbitrary} $P \in
\Pic(A)$. Thanks to (3.5.1), the equality
(3.5.2) must likewise hold for all $P$. In
other words:
$$H^i(A, \O_A(\Theta) \otimes P \otimes \J) = 0 \ \
\text{ for all }
\ i \ge 0 \ \text{ and all } \ P \in \Pic(A). $$ But it
follows from Mukai's theory \cite{Muk} of the Fourier
functor on an abelian variety that if $F$ is a coherent
sheaf on $A$ such that $H^i(A, F \otimes P) = 0$ for all
$i \ge 0$ and all $P \in
\Pic(A)$, then $F = 0$. Thus we have a contradiction.
\qed \enddemo

\vskip 15pt
\ni 

\vskip 10pt
\bl
\settabs\+University of Illinois at Chicago and stil extra
&University of calirofnia Los now is the time
\cr
\+ Lawrence EIN &Robert LAZARSFELD\cr
\+ Department of Mathematics &Department of Mathematics \cr
\+ University of Illinois at Chicago &University of California,
Los Angeles \cr
\+ Chicago, IL  60680  &Los Angeles, CA  90024 \cr
\+ e-mail: U22425\@math.uic.edu  &e-mail: rkl\@math.ucla.edu \cr

\end
{\bf References}
\Refs \nofrills{ }
\parskip 5pt

\widestnumber\key{KMori}

\bl





\ref
\key AM
\by A. Andreotti and A. Mayer
\paper On period relations for abelian integrals on
algebraic curves
\jour Ann. Scuola Norm. Sup. Pisa
\vol 21
\yr 1967
\pages 189-238
\endref

\ref
\key AD
\by E. Arbarello and C. DeConcini
\paper Another proof of a conjecture of S. P. Novikov on
 periods of abelian integrals on Riemann surfaces
\jour Duke Math. J.
\vol 54
\yr 1987
\pages 163-178
\endref

\ref
\key De
\by J.-P. Demailly
\paper $L^2$ vanishing theorems for positive line bundles
and adjunction theory
\toappear
\endref

\ref
\key GL1
\by M. Green and R. Lazarsfeld
\paper Deformation theory, generic vanishing theorems, and
some conjectures of Enriques, Catanese and Beauville
\jour Invent. Math.
\yr 1987
\vol 90
\pages 389-407
\endref

\ref
\key GL2
\bysame
\paper Higher obstructions to deforming cohomology groups of
line bundles
\jour Jour. of A.M.S.
\yr 1991
\vol 4
\pages 87-103
\endref

\ref
\key Ka
\by Y. Kawamata
\paper Characterization of abelian varieties
\jour Comp. Math.
\vol 43
\yr 1981
\pages 253-276
\endref

\ref
\key KMM
\by Y. Kawamata, K. Matsuda and K. Matsuki
\paper Introduction to the minimal model program
\inbook Algebraic Geomtery, Sendai
\jour Adv. Study in Pure Math.
\vol 10
\yr 1987
\publ Kinokuniya-North Holland
\pages 283-360
\endref


\ref
\key KV
\by Y. Kawamata and E. Viehweg
\paper On a characterization of abelian varieties
in the classification theory of algebraic varieties
\jour Comp. Math.
\yr 1981
\vol 41
\pages 355-360
\endref

\ref
\key K
\by G. Kempf
\paper On the geometry of a theorem of Riemann
\jour Ann. of Math.
\vol 98
\yr 1973
\pages  178-185
\endref

\ref
\key Kol1
\by J. Koll\'ar
\paper Higher direct images of dualizing sheaves, I
\jour Ann. of Math
\yr 1986
\vol 123
\pages 11-42
\endref

\ref
\key Kol2
\bysame
\paper Shaferevich maps and plurigenera of algebraic
varieties
\jour Invent. Math.
\yr 1993
\vol 113
\pages 177-215
\endref

\ref
\key Kol3
\bysame
\book Shafarevich Maps and Automorphic Forms
\yr 1995
\publ Princeton Univ. Press
\endref

\ref
\key Kol4
\bysame
\paper Singularities of pairs
\toappear
\endref

\ref
\key LB
\by H. Lange and C. Birkenhacke
\book Complex Abelian Varieties
\publ Springer
\yr 1992
\endref

\ref
\key Mori
\by S. Mori
\paper Classification of higher dimensional varieties
\inbook Algebebraic Geometry, Bowdoin
\jour Proc. Symp. Pure Math.
\vol 46
\yr 1987
\pages 269-332
\endref

\ref
\key Muk
\by S. Mukai
\paper Duality between $D(X)$ and $D(\hat X)$, with
application to Picard sheaves
\jour Nagoya Math. J.
\yr 1981
\vol 81
\pages 153-175
\endref

\ref
\key Nak
\by M. Nakamaye
\paper Reducibility of principaly polarized
abelian varieties with highly singular polarization
\toappear
\endref

\ref
\key SV
\by R. Smith and R. Varley
\paper Multiplicity $g$ points on theta diviors
\toappear
\endref

\ref
\key Z
\by Qi Zhang
\paper Golbal holomorphic one-forms and
Euler characteristics on projective manifolds
with ample canonical bundles
\toappear
\endref


\endRefs
